\documentclass[journal]{IEEEtran}
\usepackage[pdftex]{graphicx}
\usepackage{epstopdf}
\usepackage{caption}
\usepackage{subcaption}
\usepackage{subfig}
\usepackage{cite}
\DeclareGraphicsExtensions{.pdf,.jpeg,.png}
\usepackage{psfrag,amsfonts, amssymb}
\usepackage[cmex10]{amsmath}
\usepackage{algorithmic, enumerate}
\usepackage{cases}
\usepackage{url}

\usepackage{tabularx}

\begin{document}
%
\title{Path-based Iterative Reconstruction (PBIR) for X-ray Computed Tomography}
%
%
%

\author{Meng~Wu, 
        Andreas~Maier,
        Qiao~Yang,
        and~Rebecca~Fahrig
\thanks{M. Wu and R. Fahrig are with the Department
of Radiology, Stanford University, USA e-mail:mengwu@stanford.edu.}
\thanks{A. Maier and Q. Yang are with Pattern Recognition Lab, Friedrich-Alexander University of Erlangen-Nuremberg, Germany.}
}

\maketitle

\begin{abstract}
Model-based iterative reconstruction (MBIR) techniques have demonstrated many advantages in X-ray CT image reconstruction. The MBIR approach is often modeled as a convex optimization problem including a data fitting function and a penalty function. The tuning parameter value that regulates the strength of the penalty function is critical for achieving good reconstruction results but difficult to choose. In this work, we describe two path seeking algorithms that are capable of efficiently generating a series of MBIR images with different strengths of the penalty function. The root-mean-squared-differences of the proposed path seeking algorithms are below 4 HU throughout the entire reconstruction path. With the efficient path seeking algorithm, we suggest a path-based iterative reconstruction (PBIR) to obtain complete information from the scanned data and reconstruction model.  
\end{abstract}


%
\IEEEpeerreviewmaketitle

\section{Introduction}

The model-based iterative reconstruction (MBIR) methods for 3D computed tomography (CT) offer numerous advantages such as the potential for improved image quality and reduced dose, as compared to the conventional filtered back-projection (FBP) method \cite{Thibault2007a, Tang2009}. The statistical iterative reconstruction problem may be formulated in the Bayesian framework as a maximum a posteriori (MAP) or maximum likelihood (ML) estimation. Assuming the X-rays are monochromatic, and the detector photon counts follow the Poisson distribution, the maximum likelihood (ML) based parameter estimation is well understood as the linear Poisson regression problem \cite{Fessler:2000sl, Elbakri:2002qa, Nuyts2013}. However, the prior distribution of the object is unknown in most cases. In practice, the straightforward solution to the maximum likelihood problem generates noisy and incorrect reconstructions, especially when the number of projections is small, or the data is very noisy. A common remedy is to add a restriction/regularization to the ML solutions as an approximation of the prior distribution of the object. The reconstruction problem is then formulated as a penalized maximum likelihood problem (PML)
\begin{equation}
\begin{aligned}
\mu & = \underset{\mu \ge 0}{\text{argmax}} \Psi( \mu ) - \beta h ( \mu ), \\
\end{aligned}
\label{eqn:pml}
\end{equation}
where $\Psi( \mu )$ is the log-likelihood function, $h( \mu )$ is the penalty function (also known as regularization), $\mu$ is the reconstruction, and $\beta$ is the tuning parameter that regulates the strength of the penalty function. 

Over the last twenty years, there has been extensive research on finding the most suitable penalty function for CT reconstructions. The penalty function can be chosen to enforce the image smoothness and sparsity \cite{Elbakri:2002qa, Sidky2008, Pfister2014}, or designed to produce spatially invariant or uniform point spread function and noise \cite{Fessler:1996bs, Stayman2000, Cho2015}. Prior CT volumes have also been considered for regularization of the iterative reconstruction in extremely undersampled view cases \cite{Tang2009}. Some recent work also suggests modification of the maximum likelihood function to generate images with certain properties \cite{chang:2015}. 

Unfortunately, not only is the perfect penalty function unknown, but also the right strength of the penalty function is difficult to select for different cases. In the PML reconstruction, different values of tuning parameter $\beta$ generate different reconstructed images (solutions to the slightly different optimization problems). In fact, the values of the tuning parameter $\beta$ ($0 \le \beta \le \infty $) in Eqn. \eqref{eqn:pml} produce a series of reconstructions $\mu(\beta)$. It is well known that the value of the tuning parameter is critical to the reconstruction results \cite{Tang2009, Wang2014a, Han2015}. For example, if $\beta$ is too small, the regularization is not strong enough to suppress noise and artifacts; if $\beta$ is too big, the image is over blurred and even exhibits patchy behavior. 

To the best of our knowledge, there is no perfect way to choose the value of $\beta$ that would lead to the best reconstruction with maximum clinical utility \cite{Han2015}. Trial-and-error or exhaustive search of the best tuning parameter are often seen in the literature \cite{Wang2014a, Stayman2013, Wu2014a}. Approximation methods based on the analytical solution of the MBIR have also been studied \cite{Dang2014b}. On the other hand, the MBIR requires solving a very large scale optimization problem. Exploration of fast optimization solvers such as accelerated first order methods and alternating direction method of multipliers (ADMM) has been one of the hottest topics in recent years \cite{Niu:2012fk, Kim2013a, Nien2014}. Although 50 to 100 iterations using the state-of-art optimization solvers can produce an accurate enough solution for one tuning parameter value, directly computing multiple solutions ($\mu(\beta)$) via numerical optimization would not be suitable in practice.

\begin{figure}[!t] 
   \centering
   \includegraphics[width=0.45\textwidth]{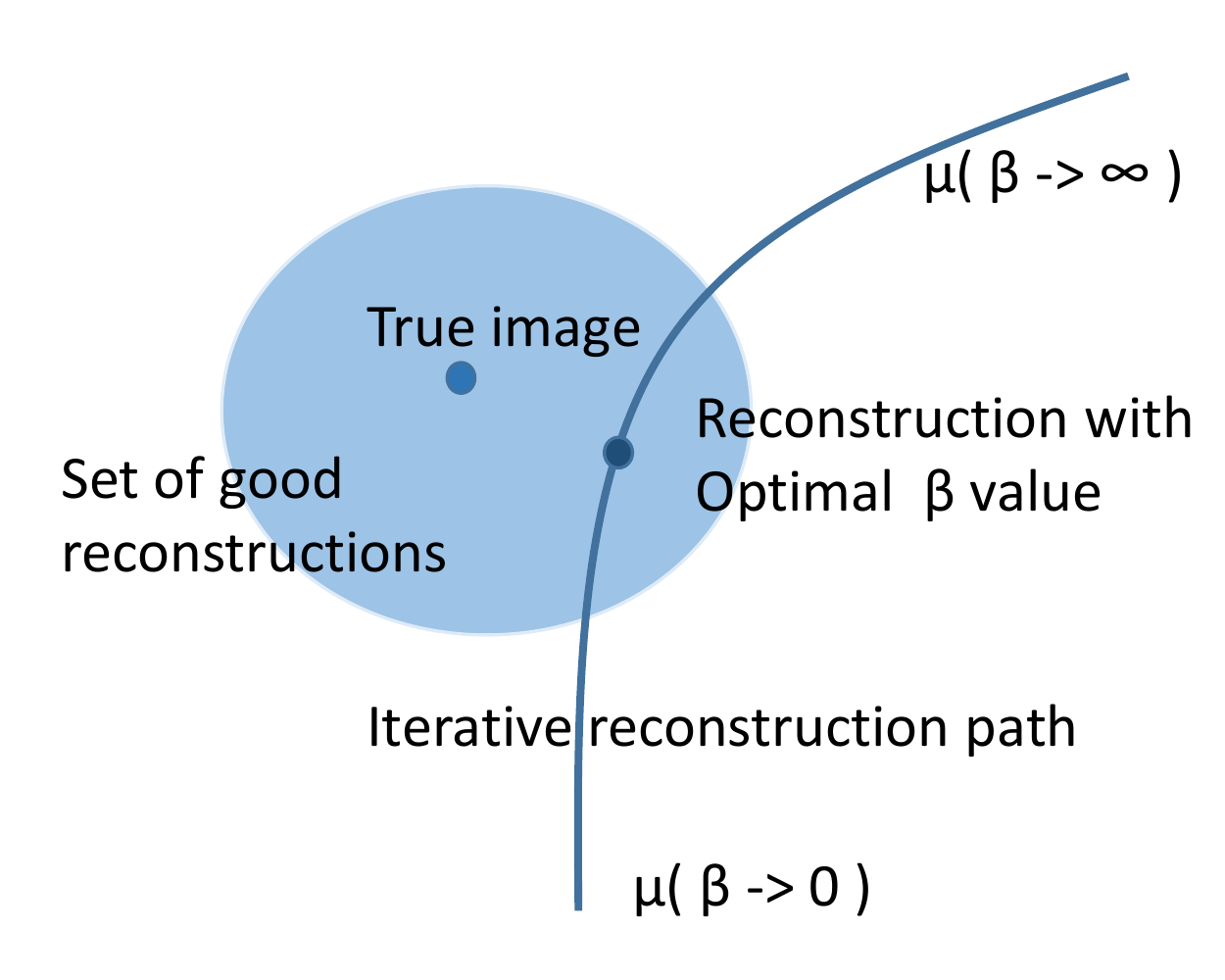} 
   \caption{ Reconstruction path of the penalized maximum likelihood method. }
   \label{fig:rp}
\end{figure}

Instead of focusing on finding the single optimal tuning parameters, we suggest evaluating the reconstruction path of the MBIR. Figure \ref{fig:rp} shows the relationship between the true image and the PML reconstructions with different tuning parameter values. The true image is the desired reconstruction result that provides the best representation of the scanned object. Due to the limitations of physics, system defects, and noise, the true image will almost always be unreachable. But around the true image, there is a set of good reconstructions, which provide sufficient information for clinical use. A good reconstruction model will have a segment of the reconstruction path within the set of good reconstructions. To get the closest iterative reconstructions to the true image requires selecting the optimal tuning parameter value prior to the reconstruction. Therefore, a method that provides the entire reconstruction path can provide complete information for a given reconstruction model and removes the burden of selection of good tuning parameters. We call the new reconstruction strategy path-based iterative reconstruction (PBIR). In the paper, we describe two path seeking algorithms to extend MBIR to PBIR.

Recently, path seeking algorithms have been studied for regularized regression problems and have shown advantages in selecting the tuning parameter\cite{Hastie2009, Friedman2012}. Many statistics packages have solutions for regression problems from the path perspective \cite{fiedman2009}. Friedman proposed a fast generalized path seeking (GPS) algorithm that produces solutions closely approximating the path of the constrained regression problems for non-identical form and penalty function between $L_1$ and $L_2$ norms \cite{Friedman2012, Friedman2007pc}. The generalized path seeking algorithm uses the ratio of the gradients to update one regression variable each time, which is accurate and suitable for a small regression problem. The dual path algorithm proposed by Tibshirani and Taylor operates in a single, unified framework that allows the $L_1$-based regularization to be completely general \cite{Tibshirani2011, Arnold2014, zhou2011zl}. 

The path seeking in the image reconstruction problem is very similar to the regularized regression problem. The two main differences between the image reconstruction and regularized regression problems are: 1. the scale of the path seeking in image reconstruction is much greater than the regularized regression problem.  2. in the regularized regression, the tuning parameter $\beta$ is often set to a very large value such that the variables are equal to zero. The end of the path will have the tuning parameter equal to zero, which leads to the ordinary least-square problem. In image reconstruction, the range of the $\beta$ values that make useful images is much smaller. In this paper, we describe two path seeking algorithms for the large-scale image reconstruction problem. Evaluations and comparison of the path seeking algorithms in the image domain and via noise power spectrum are presented. We introduce the path-based iterative reconstruction and discuss its potential benefits. 

\section{Methods}

In this study, we consider the penalized weighted least-squares (PWLS) algorithm \cite{Elbakri:2002qa}
\begin{equation}
\centering
\begin{aligned}
\mu & = \underset{\mu \ge 0}{\text{argmin}} \frac{1}{2} \sum_i w_i ( [ P\mu]_i - l_i )^2   + \beta h( \mu )
\label{eqn:pwlsobj}
\end{aligned}
\end{equation}
where $P$ denotes the system matrix for the data acquisition geometry, $l_i$ denotes the logged normalized projection of the $i$th ray, and $w$ is the least-squares weight to account for the noise level in the X-ray projection data. In this paper, we used the penalized least-square notation for simplicity
\begin{equation}
\begin{aligned}
 \text{ minimize} & \quad \frac{1}{2} \| A \mu - y \|_2^2 + \beta h ( \mu ) =  g(\mu) + \beta h ( \mu ) \\
 \text{subject to} & \quad \mu \ge 0,
 \end{aligned}f
 \label{eqn:pls}
\end{equation} 
where $A = W^{1/2} P $, $ y = W^{1/2} l $, and $g(\mu)$ denotes the least-squares part. $W$ is the diagonal matrix containing $w_i$.

In this paper, we present two accurate and efficient algorithms to compute the reconstruction path of the MBIR. The first approach is based on the generalized path seeking algorithm using the ratio-of-gradients information \cite{Friedman2012, wu2015ps, Wu:2015tps}. The second approach uses the direction-of-gradient constrained optimization to obtain proper path seeking direction and step size.

\subsection{Ratio-of-Gradients Search}

The generalized path seeking algorithm uses the ratio-of-gradients to update one of the regression variables with a fixed step size \cite{Friedman2012}. The basic idea of selectively updating pixel values in our first approach is the same as generalized path seeking algorithm. The negative gradient of the least-squares term for the penalized least-squares problem from Eqn. \eqref{eqn:pls} is 
 \begin{equation}
\nabla_j g(\mu) = - \left[ A^T  ( A \mu - y ) \right]_j,
\end{equation}
and the ratio-of-gradients is 
 \begin{equation}
\lambda_j = \frac{ \nabla_j  h (\mu) }{ | \nabla_j  g(\mu) | }.
\label{eqn:rog}
\end{equation}
The magnitude of $\lambda_j$ reflects the relative strengths of two functions in Eqn. \eqref{eqn:pls} for the $j$th pixel. For example, large $|\lambda_j|$ means the penalty function has a stronger effect than the least-squares function on the $j$th pixel. When the tuning parameter $\beta$ increases, the $j$th pixel is more likely to be changed than the other pixels with smaller $|\lambda|$. The sign of $\lambda$ indicates the direction of the change. In the path seeking update, the algorithm can update the fraction (10\% - 20 \%) of pixels with the largest absolute values of the ratio-of-gradients along the direction of $\lambda$ by a fixed amount, chosen to be 1 - 2 HU here \cite{wu2015ps}. 

The generalized path seeking algorithm has been shown to be very accurate with a single variable update and small step size. However, using only the ratio-of-gradients updates for the MBIR path seeking is not accurate, because the fixed size update is clearly non-optimal, and errors will accumulate as the path images move away from the initial image. We introduce two techniques to improve the efficiency and accuracy of the ratio-of-gradients based path seeking called target direction and intermediate optimization. 

As discussed in the previous section, the path seeking for the image reconstruction problem does not have to run the tuning parameter from zero to infinity. There is a rough range of the tuning parameter value that contains good reconstruction results. Choosing the range is much easier than selecting a single optimal value. We can first do two iterative reconstructions and use them as the start and end points of the path seeking. For example, if we want to do path seeking from  $\mu( \beta_1 )$ to $\mu( \beta_2 )$, the target direction of the $j$th pixel is then defined as 
 \begin{equation}
\begin{aligned}
d_j = \text{sign}\left\{ \mu_j( \beta_2 ) - \mu_j( \beta
_1 ) \right\}, 
\label{eqn:target_dir}
 \end{aligned}
\end{equation}
for all path images within the search range. When $\beta_1 $ and $ \beta_2$ are close, we assume the path seeking direction is the same as the target direction
 \begin{equation}
\begin{aligned}
 \frac{ \partial( \mu_j(\beta) )}{ \partial \beta } \cdot d_j > 0
 \end{aligned}
\end{equation}
for all $\beta \in [\beta_1, \beta_2]$. Then the pixel is updated only if its target directions\ is the same as its ratio-of-gradients direction as defined in Eqn. \eqref{eqn:rog}. This constraint is equivalent to saying the path of each pixel is locally monotonic in $\beta$. Note that, this monotonic assumption is not true in either theory or practice. But introducing the target gradient to constrain the updating pixel set can improve both path seeking efficiency and accuracy in practice. 

The second technique to improve the accuracy is adding actual optimization steps in the path seeking process, since the path images are supposed to be the solutions of a series of optimization problems, which differ only in $\beta$ values. Thus, the tuning parameter value can be estimated using the Karush-Kuhn-Tucker (KKT) conditions from the new path image \cite{Boyd2010}. The KKT condition for the penalized least-squares problem is
\begin{equation}
\begin{aligned}
\nabla_j g( \mu^{\ast} ) + \beta \cdot \nabla_j h ( \mu^{\ast} ) - \eta_j^{\ast} &  = 0  \\
\eta_j^{\ast} \cdot \mu_j^{\ast} &  = 0  & \text{for all} \quad j, \\
\eta_j^{\ast} & \ge 0 
\end{aligned}
\end{equation}
where $\mu^{ast}$ is the solution of a PWLS problem, and $\eta_j^{\ast}$ is the Lagrange multiplier of the non-negative constraint in Eqn. \eqref{eqn:pls}. From the convex optimization point of view, with a given tuning parameter, solving the KKT conditions (if strictly feasible) gives both the primary and dual solutions. From the path seeking point of view, with a prime solution of the KKT conditions, the tuning parameter value $\beta$ can be estimated by
\begin{equation}
\beta \approx \text{Median} \left\{ \frac{\nabla_j g( \mu^{\ast} )}{ \nabla_j h( \mu^{\ast} )}, \quad \forall j: \mu_j > 0  \right\}.
\label{eqn:beta}
\end{equation}
Using the median is more numerically robust than using the mean. When $\beta$ is large, the image is often very smooth (assuming the penalty function for roughness). The magnitude of $\nabla_j h( \mu^{\ast} )$ is very small which can cause unexpectedly large $\beta$ in the estimation. With a way of estimating $\beta$, we can add several optimization iterations before or after each ratio-of-gradients updating step.

\begin{table}[!t]
\caption{Pseudo code for path seeking algorithm using ratio-of-gradients (PS-ROG)}
\normalsize
\centering
\begin{tabular}{|p{0.46\textwidth}|}
\hline

Reconstruct two images $\mu( \beta_1 )$ and $\mu( \beta_2 )$ for selected path range $[\beta_1, \beta_2]$ \\

Set the initial path image $\mu = \mu( \beta_1 ) $ \\

Loop \{ 

\begin{enumerate}
\item Estimate $\beta$ using Eqn. \eqref{eqn:beta}

\item Run $\mu = \mu - \alpha \left(  \nabla g( \mu ) + \beta \cdot \nabla h( \mu ) \right)$ several times

\item Compute target direction $d_j$ using Eqn. \eqref{eqn:target_dir}

\item Compute ratio-of-gradients $\lambda_j = \nabla_j h(\mu) / | \nabla_j g(\mu) | $
 
\item Find $S = \{ j |  h_j(\mu) \cdot g_j(\mu) >  0   \}$
 
\item If ($S$ is not empty) \{

\item \quad Update $\mu_j = \mu_j + \Delta v \cdot \text{sign}( \lambda_j ), \forall  j \in S  $

\item \} Else \{

\item \quad  If ( $\lambda_j \cdot d_j < 0 $  ), then $ \lambda_j = 0 $

\item \quad  Find $t$ such that $\text{Prob}\{ | \lambda_j | \ge t \} \le p $

\item \quad  Update $\mu_j = \mu_j + \Delta v \cdot d_j, \forall | \lambda_j | \ge t $

\item \}
\end{enumerate}

\}\\

Until $\| \mu - \mu(\beta_2) \| $ stops decreasing. \\

\hline
\end{tabular}
\label{tab:ps_rog}
\end{table}

Table \ref{tab:ps_rog} presents the pseudo-code of the ratio-of-gradient based path seeking algorithm. The variable $p$ is the percentage of the pixels updated in each iteration, and $\Delta v > 0$ is a small increment value, e.g. 1 - 2 HU. Lines 1-2 are gradient descent-based optimization steps to draw $\mu$ closer to the correct reconstruction path. At each iteration, the pixel-wise ratio-of-gradients and the target direction are computed for the current path image $\mu$. Line 7 updates the pixel if both gradients are in the same direction. Line 9 ensures the algorithm only considers the rest of the pixels that have the same updating direction as the target direction. Line 10 selects the pixels that have the largest ratios of the gradients. The selected pixel is then incremented by a small fixed amount ($\Delta v$) in the target direction.

\begin{figure}[!t] 
   \centering
   \includegraphics[width=0.45\textwidth]{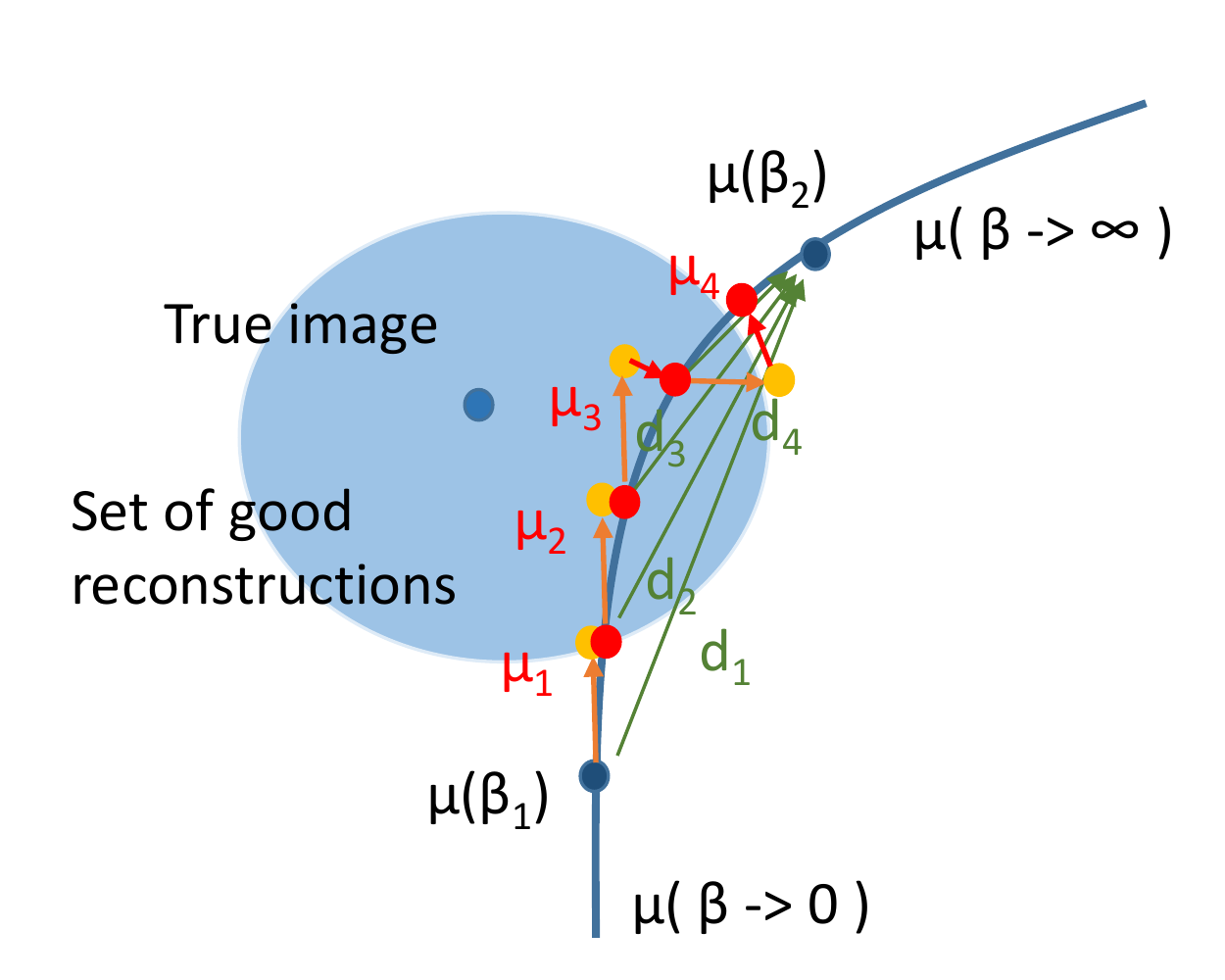} 
   \caption{ 2D illustration of the ratio-of-gradients path seeking algorithm. The red points ($\mu_1 \dots \mu_4$) are the path images computed by the true path seeking algorithm. }
   \label{fig:ps_rog}
\end{figure}

Figure \ref{fig:ps_rog} shows a 2D illustration of the ratio-of-gradients based path seeking algorithm. The green arrows $d_1 \dots d_4$ are the target directions toward final path image $\mu( \beta_2 )$; the orange arrows are the fixed size updates to the image (line 11). In each iteration, one out of two variables was updated by $\Delta_v$ in the same direction as $d$ and $\lambda$. The images from the previous path seeking step (orange points) are corrected by the minimization step (red 
arrows) to more accurate path images (red points) before the next path seeking iteration.

Note that, the ratio-of-gradients based path seeking can be done for both increasing or decreasing tuning parameter values. For decreasing $\beta$ values, we just need to swap the initialization and target images,  and use $d_j = \text{sign}(\mu_j( \beta_1 ) - \mu_j )$ and $\lambda_j = \nabla_j g(\mu) / | \nabla_j h(\mu) | $. In our experiments, we found the path seeking algorithm with increasing $\beta$ value is more stable and accurate than decreasing $\beta$. Because when $\beta$ is large, the image is very smooth and gradient of the regularization ($\nabla_j h(\mu)$) is close to zero.  

The ratio-of-gradient based approach selects a subset of pixels for a fixed step size update. However, the accuracy of the selected updating set of pixels may suffer from the previously incorrect path solution and accumulated path seeking errors. Moreover, the step size is clearly not optimized for every pixel and this may cause overstepping or understepping in the update. The disadvantages of the ratio-of-gradients method motivate us to investigate another path seeking approach that updates the entire image simultaneously. 

\subsection{Direction-of-Gradient Search }

The second approach is called the direction-of-gradient based path seeking. It uses the direction of one gradient function to constrain the optimization problem thereby encouraging the image to change in the desired direction. For example, if we want to seek the path of increasing strength of the penalty function $h( \mu )$, then we would like to encourage the optimization updates (i.e. gradient descent) to go in the same direction as the $\nabla h( \mu )$. Let us consider adding a linear inequality constraint to the penalized least-squares problem as 
\begin{equation}
\begin{aligned}
 \text{ minimize} & \quad   g( \mu )  + \beta_1 h( \mu )  \\
 \text{subject to} & \quad  \mu \ge 0  \\
 & \quad ( \mu_j - \hat{\mu}_j ) \cdot \nabla_j h( \hat{\mu} ) \le 0 \quad \forall j,
 \end{aligned}
 \label{eqn:opt_inactive_lc}
\end{equation}
where 
\begin{equation}
\hat{\mu} = \underset{\mu \ge 0}{\text{argmin}} \quad g( \mu )  + \beta_1 h( \mu ).
\label{eqn:mu_hat}
\end{equation}
The second linear inequality constraint in Eqn. \eqref{eqn:opt_inactive_lc} is inactive because the $\hat{\mu}$ is already optimal for the nonnegative constrained penalized least-squares problem.

If we slightly increase $\beta_1$ to $\beta_2$ in the direction-of-gradient constrained penalized least-squares problem \eqref{eqn:opt_inactive_lc} as
\begin{equation}
\begin{aligned}
 \text{ minimize} & \quad  g( \mu )  + \beta_2 h( \mu )  \\
 \text{subject to} & \quad  \mu \ge 0  \\
 & \quad ( \mu_j - \hat{\mu}_j ) \cdot \nabla_j h( \hat{\mu} ) \le 0 \quad \forall j
 \end{aligned}
 \label{eqn:opt_active_lc}
\end{equation}
and keep the $\hat{\mu}$ same as in Eqn. \eqref{eqn:mu_hat}, the new solution will be suboptimal for the penalized least-squares problem \eqref{eqn:pls} with $\beta_2$. But the solution of the problem \eqref{eqn:opt_inactive_lc} is still close to the solution of the reconstruction problem because increasing the strength of $h(\mu)$ and the direction-of-gradient constraint have very similar effects. To solve the direction-of-gradient constrained problem, we can simply apply a projection onto convex sets (POCS) step \cite{Sidky2008} in the optimization step (e.g. gradient descent) as described in Table \ref{tab:pls_pocs}. The POCS step will encourage updates of the image that favor minimizing $h(\mu)$, which increases the path seeking efficiency within the optimization framework.

\begin{table}[!t]
\caption{Projection onto convex sets solver for problem \eqref{eqn:opt_active_lc} }
\normalsize
\centering
\begin{tabular}{|p{0.45\textwidth}|}
\hline
Warm start $\mu = \hat{\mu}  $ \\
Loop \{ 
\begin{enumerate}
\item $\mu = \mu - \alpha \left(  \nabla g( \mu ) + \beta_2 \cdot \nabla h( \mu ) \right)$ 
\item Set $\mu_j = 0, \quad  \forall \quad \mu_j \le 0 $
\item Set $\mu_j =  \hat{\mu}_j, \quad  \forall  \quad ( \mu_j - \hat{\mu} _j ) \cdot \nabla_j h( \hat{\mu} ) \ge 0$.
\end{enumerate}
\}\\
\hline
\end{tabular}
\label{tab:pls_pocs}
\end{table}

For every iteration in Table \ref{tab:pls_pocs}, we need to calculate both $\nabla g(\mu)$ and $\nabla h(\mu)$ in the first gradient descent step. The computational cost of $\nabla g(\mu)$ is very high because it requires forward and backward projections. It is not computationally efficient to spend too much computing power on the direction-of-gradient constrained problem, which only gives a fast approximation to the path image at each new $\beta$. Therefore, we propose to split the variable such that a surrogate function of $g(\mu)$ with a simple gradient formula can be used for the direction-of-gradient constrained sub-problem. The penalized least square problem can be modified by variable splitting as in \cite{ Nien2014}:
\begin{equation} 
\begin{aligned}
&\text{minimize} \quad  \frac{1}{2} \| z - y \|^2_2 + \beta h ( \mu ) \\
&\text{subject to} \quad  A \mu = z 
\label{eq:pls_split}
\end{aligned}
\end{equation} 
We further use the Augmented Largrangian (AL) to solve the constrained minimization problem:   
\begin{equation}
\mathcal{L}_A( \mu, z, t; \rho) \triangleq \frac{1}{2} \| z - y \|^2_2 + \beta h( \mu ) + \frac{\rho}{2} \| A \mu - z - t \| ^2_2
\end{equation}
where $t$ and $\rho > 0$ are the corresponding AL penalty parameters. The alternating direction method of multiplier (ADMM) solver can be summarized as \cite{Chamboll2011, Nien2014}:
\begin{equation}
\begin{cases}
&s^{(k+1)} = \rho A^T (A \mu^{(k)} - y) + (1 - \rho) v^{(k)}\\
&\mu^{(k+1)} = \underset{\mu \ge 0}{\operatorname{argmin}}\{\beta h(\mu) + \frac{\rho}{2t} \| \mu - \mu^{(k)} + s^{(k+1)} \| ^2_2\} \\
&v^{(k+1)} = \frac{\rho}{\rho+1}A^T (A \mu^{(k+1)} - y) + \frac{1}{\rho+1} v^{(k)},
\end{cases}
\label{eqn:admm}
\end{equation}
where $v \triangleq A^T(z - y)$ is the backprojection of the split residual. In the second step, the penalty function $h( \mu )$ is only interacting with a simple quadratic form of $\mu$. We modify the second step with additional direction-of-gradient constraint as 
\begin{equation}
\begin{aligned}
 \text{ minimize} & \quad  \beta h( \mu ) + \frac{\rho}{2t} \| \mu - \mu^{(k)} + s^{(k+1)} \| ^2_2 \\
 \text{subject to} & \quad \mu \ge 0 \\
 & \quad  ( \mu_j - \mu_j^{(k)} ) \cdot \nabla_j h( \mu^{(k)}) \le 0 \quad \forall j.
 \end{aligned}
 \label{eqn:admm_dog}
\end{equation}
To solve this sub-problem, one can use simple gradient descents and POCS similar to Table \ref{tab:pls_pocs}. In this way, we can run multiple iterations to solve the direction-of-gradient constrained optimization problem \eqref{eqn:admm_dog}. The strength of $h( \mu )$ is increased under the optimization framework while the computational load of solving the sub-problem is still reasonable. As for the ratio-of-gradients based approach, additional ordinary ADMM steps in Eqn. \eqref{eqn:admm_dog} can be used to improve the accuracy of the path seeking. The direction-of-gradient based path seeking algorithm is summarized in Table \ref{tab:ps_dog}.

\begin{table}[!t]
\caption{Pseudo code for path seeking algorithm using direction-of-gradient (PS-DOG)}
\normalsize
\centering
\begin{tabular}{|p{0.45\textwidth}|}
\hline
Reconstruct an image $x = x_{ \beta_1 } $ \\
Loop \{ \
\begin{enumerate}
\item Execute modified ADMM steps in Eqn. \eqref{eqn:admm} and Eqn. \eqref{eqn:admm_dog}.
\item Execute normal ADMM steps in Eqn. \eqref{eqn:admm} for several iterations.
\item If $\| x - x_{\beta_1} \| $ is not increasing, the increase $\beta$.
\end{enumerate}
\}\\
Until $\beta = \beta_2$. \\
\hline
\end{tabular}
\label{tab:ps_dog}
\end{table}

\subsection{Ordered Subset Acceleration }

The ordered subsets method is commonly used to accelerate the speed of convergence of iterative reconstruction solvers \cite{hudson1994, Elbakri:2002qa, Nien2014}. The ordered subsets method uses a small fraction of the projection data to estimate the data fitting function and its gradient in the optimization step. The path seeking methods also need to compute the gradient of the data fitting function in each iteration. The computational load of repeated projection and backprojection of full projection data can be reduced by only using a subset of the projections that are equally sampled in projection angles. 

Common choices of the number of ordered subsets for a clinical CT system is between 20 and 40 (20-50 projections per subset) depending on the scan geometry and optimization solver \cite{ Elbakri:2002qa, Nien2014, Kim2015}. In general, the number of projections per subset in the path seeking method must be larger than the direct optimization in order to have more accurate estimation of the gradient. For standard optimization problems, the image is assumed to converge to a single point. In contrast, the path seeking algorithms encourage the image to step away from the current convergence point to a new convergence point. The path seeking algorithms are therefore naturally less stable than the optimization algorithms.

For the ratio-of-gradients based path seeking algorithm, the suitable number of ordered subsets are 5 - 10, because the ordered subset errors in the ratio-of-gradients will have the accumulated path seeking errors. Updating pixels in the wrong order may cause the path seeking solution to diverge from the correct path. The optimization steps in the ratio-of-gradients path seeking method can have the number of ordered subsets as large as in the direct optimization method \cite{Elbakri:2002qa}.  

For the direction-of-gradient based method, the path seeking is under the framework of the constrained optimization problem that is more robust to ordered subset errors than the fixed step size update. Also, the direction-of-gradient updates use the quadratic surrogate function, which contains information of the data fitting errors ($s$ in Eqn. \eqref{eqn:admm_dog}). Thus, the suitable number of ordered subsets for the direction-of-gradient path seeking method is between 10 and 20. In order to execute alternatively between the normal and modified ADMM optimization steps, the additional optimization steps need to have the same number of ordered subsets\cite{Nien2014}. 

\section{Simulations}

A 64-slice clinical diagnostic CT scanner geometry (LightSpeed, GE Healthcare, Waukesha, WI) was used in the simulations. There are 984 projections per 360 degrees circular rotation, and the detector size is  888$\times$64. The reconstructed image size 512$\times$512$\times$30 with in-plane spacing of 0.7 - 0.9 mm$^2$ depending on the phantom size. The slice thickness is 1 mm for all the reconstructions. 

Three numerical phantoms were used in this work: two XCAT phantoms of abdomen and thorax and a water cylinder phantom with a diameter of 32 cm (body CTDI phantom). The voxel size of the XCAT phantom is $0.6$ mm isotropic, and the voxel size of the water cylinder phantom is $0.4$ mm isotropic. All projection data were simulated in an axial scanning mode using a 120 kVp polychromatic spectrum. Simulated projections of the XCAT phantom were generated assuming an exposure of approximately 100 mAs, corrsponding to 2$\times 10^5$ photons per unattenuated ray. The water cylinder projections were simulated at 50 mAs and 100 mAs.

The simulated projection data are reconstructed using both filter-backprojection (FBP) and penalized weighted least-squares (PWLS) methods. We used the convex edge-preserving Huber function as the penalty function for image roughness. The transition value from quadratic to linear regions is set to 5 Hounsfield units (HU), which has been reported to provide a good trade-off between soft-tissue contrast and noise reduction \cite{Wang2013}. The direct optimization solutions of the PWLS reconstruction were achieved using the 20 ordered-subsets linearized augmented Lagrangian method with 50 iterations \cite{Nien2014}. 

\begin{figure*}[!t]
\centering
\includegraphics[width=\textwidth]{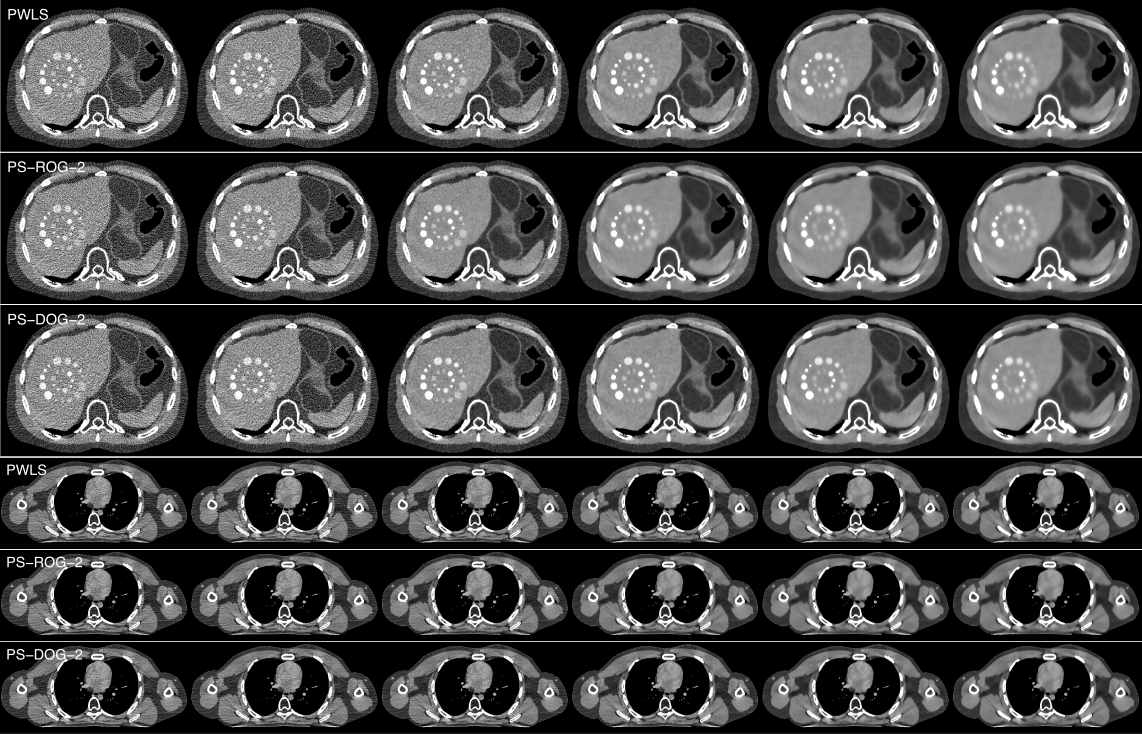}
\caption{ Six path images in the PWLS reconstructions path computed using the direct optimization (top), the proposed PS-ROG-2 algorithm (middle), and the proposed PS-DOG-2 algorithm (bottom) using the abdomen and the thorax XCAT phantoms. The tuning parameters of the path images are evenly log-spaced from $5\times10^3$ to $2\times10^5$. The display window is [0 150] HU. }
\label{fig:recon_path}
\end{figure*}

The proposed path seeking methods were used to generate path images of the PWLS reconstruction with $\beta$ values from $5\times 10^3$ to $2 \times 10^5$ for all cases. The range of tuning parameters produces reconstructions ranging from very noisy to over smoothed. The optimal choices of tuning parameter that balance the trade-offs between noise reduction and resolution are within this range. Note that, the path seeking algorithm is designed to permit efficient calculation of reconstructions for monotonically changing tuning parameter, but can not be used to investigate the impact of changing parameters in the penalty function itself. To investigate the effect of different penalty function parameters, one would run independent reconstructions and their corresponding path seeking algorithms. In this validation, we only vary the path seeking tuning parameter value from small to large, which has been demonstrated to be more accurate than the opposite direction \cite{wu2015ps}. 

A total of 40 path images with roughly equal mean-absolute-differences were computed using the ratio-of-gradients path seeking method. The update percentage $p$ was set to 20\%, with step size $\Delta v $ of 1 HU. Zero to two sub-iterations of the separable quadratic surrogate (SQS) gradient descent optimization steps are executed to improve accuracy before storing each path image \cite{Elbakri:2002qa, Fessler:2000sl}. The reconstructions are denoted as PS-ROG-N, where N is the number of SQS sub-iterations per frame. The number of the ordered subsets is 5 for the path seeking steps and 20 for the gradient descent optimization steps. The number of iterations required to generate the entire sequence of path images is 50 $\times$ 2 (two initial reconstructions) + 40 (path seeking updates) + 40 $\times$ N (optimizations). The number of iterations for path seeking can be reduced by selecting a smaller number of path frames. The iterations in the reconstruction and the path seeking have similar complexity. 

For the direction-of-gradient path seeking method, a total of 40 path images with log-spaced tuning parameter values are computed and stored. The increment ratio of the tuning parameter is 1.45. Zero to two normal ADMM optimization iterations are executed between the modified direction-of-gradient constrained ADMM steps to improve the accuracy. The reconstructions are denoted as PS-DOG-N, where N is the number of  intermediate ADMM optimization iterations. The number of the ordered subsets was 10 for both modified and normal ADMM optimization steps. The number of iterations required to generate the entire sequence of path images is 50 (one initial reconstruction) + 40 (path seeking updates) + 40 $\times$ N (optimizations).

\section{Results}

\begin{figure*}[!t]
\centering
\includegraphics[width=\textwidth]{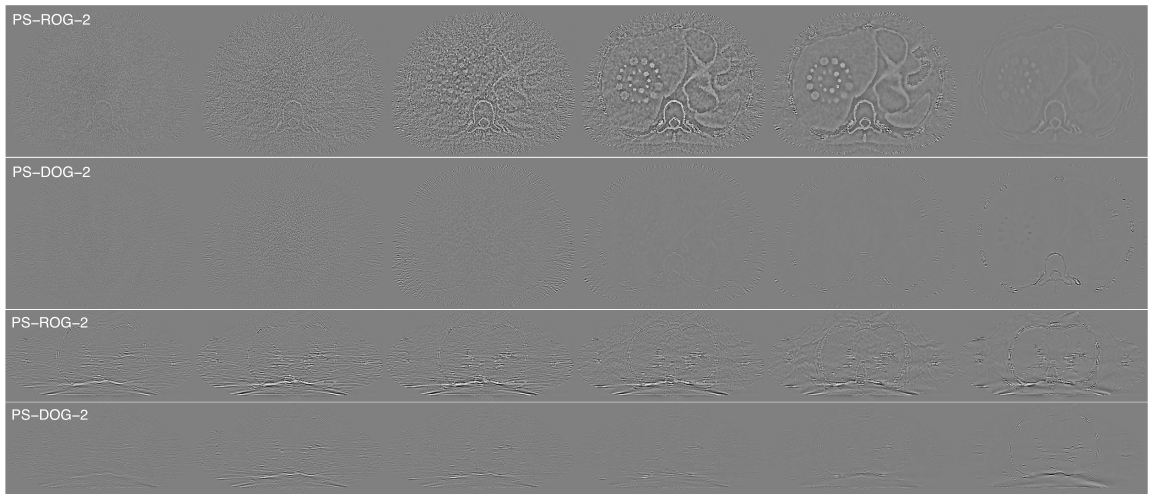}
\caption{ Difference images between direct reconstructions in Figure \ref{fig:recon_path} and the closest path seeking images. The display window is [-20 20] HU. }
\label{fig:recon_diff}
\end{figure*}

Figure \ref{fig:recon_path} shows six frames in the reconstruction path of the PWLS method using the direct optimization (ground truth) and the proposed the path seeking methods. All of the approaches provided sequences of images from noisy reconstructions to over blurred images. Both of the proposed path seeking approaches generated images that are visually similar to the path images using the direct optimization. The  soft tissue noise texture using the PS-DOG-2 method is more similar to the ground truth. Figure \ref{fig:recon_diff} shows the difference images between direct reconstructions in Figure \ref{fig:recon_path} and the closest path seeking images. The errors in the PS-ROG images are larger than in the PS-DOG images. The path seeking errors using the PS-ROG method first accumulate and then decrease because the path is also constrained by the target image. The PS-DOG method only has small errors around the edges.

\begin{figure*}[!t]
\centering
\includegraphics[width=\textwidth]{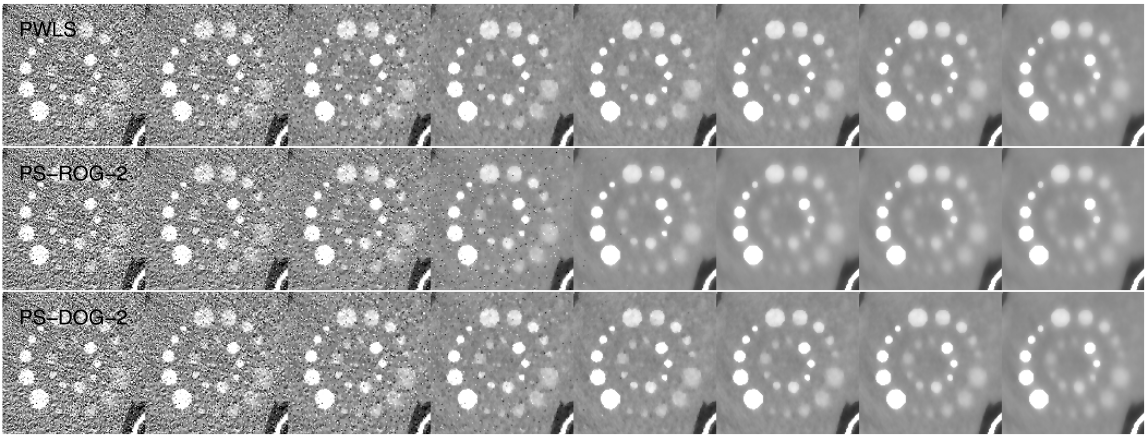}
\caption{ Eight 9 cm $\times$ 9 cm regions of interest in the PWLS reconstructions path computed using the direct optimization (top), the proposed PS-ROG algorithm (middle), and the proposed PS-DOG algorithm (bottom) .}
\label{fig:rois}
\end{figure*}

\begin{figure*}[!t]
\centering
\includegraphics[width=0.8\textwidth]{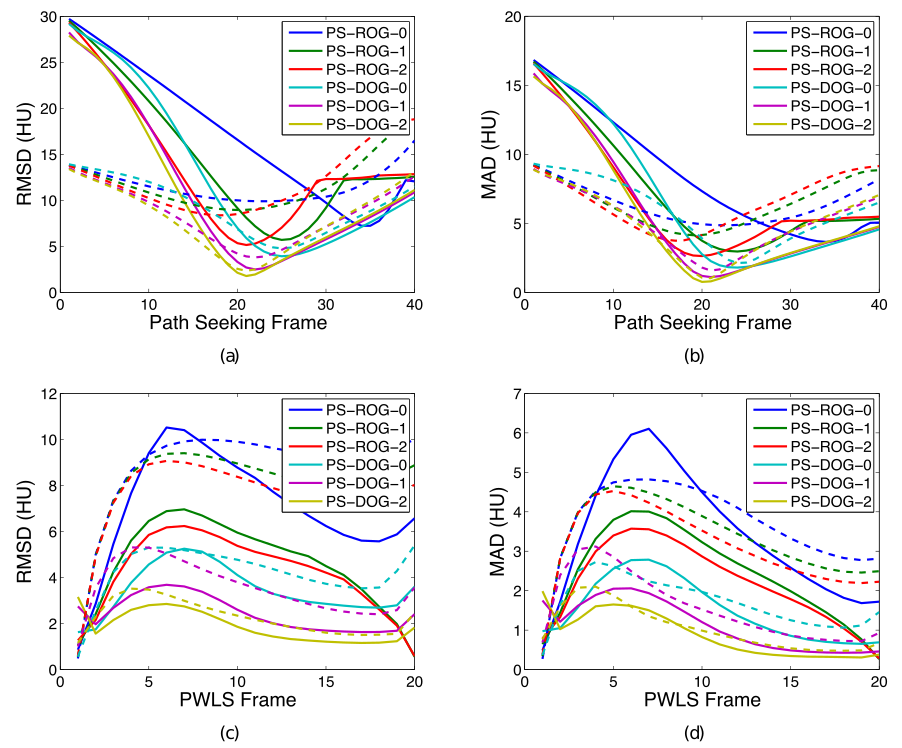}
\caption{The RMSD and MAD measurements of the entire path images generated by the proposed method compared to a directly solved PWLS images. The solid lines are the measurements of abdomen phantom, the dashed lines are the measurements of the thorax phantom. }
\label{fig:measure}
\end{figure*}

Figure \ref{fig:rois} shows 9 cm $\times$ 9 cm region-of-interest (ROI) of the path images using direct optimization and the proposed path seeking methods. The path images generated by the path seeking algorithms are similar to the direct optimization results. There are some salt and pepper noise in the 4th image using the PS-ROG method.  Those pixels are updated slower than the rest of the image because the accumulated errors in the PS-ROG method. But with the additional optimization, the errors have been corrected quickly. The PS-DOG path images are visually closer to the directly solved PWLS images than the PS-ROG, especially in the background variation.

We used the root-mean-squared-difference (RMSD) and mean-absolute-difference (MAD) as quantitive measures of the path seeking accuracy. The RMSD and MAD between the first and the last path images are 37.6/22.5 HU and 21.3/13.1 HU (RMSD/MAD) for the abdomen and thorax cases, respectively. Each path image solved by the direct optimization method is compared to the entire reconstruction paths generated by the proposed path seeking method. Figure \ref{fig:measure} (a) and (b) show examples of RMSD and MAD of one directly solved PWLS image with $\beta = 6\times10^4$ (20th frame) compared to the 40 path images. The closest path seeking image in the middle of the path has RMSD and MAD around 2-3 HU using the PS-DOG-2 method. The minimum differences using the PS-ROG method are larger than those using the PS-DOG method. The frame numbers of the closest path images are not constant when using a different number of additional optimization steps in the ratio-of-gradient based method. The frame numbers of the closest path images using the PS-DOG method is more stable and predictable because the path frames are controlled by the tuning parameter values. The estimation of the tuning parameter is not robust when the image is smooth, and therefore directly linking the PS-ROG method with the actual tuning parameter values of the reconstructed images is more difficult.

Figure \ref{fig:measure} (c) and (d) show the minimum RMSD and MAD for all of the directly solved PWLS images compared with the path seeking images. Again, the PS-DOG method is more accurate than the PS-ROG method. The worst RMSDs of the entire PS-DOG-2 path are less than 4 HU. The errors in the first half of the reconstruction path are larger than in the second half, because the tuning parameter changes in the first part of the path are mainly suppressing the noise, which is more difficult to track than the blurring effects in the later part of the path. 

\begin{figure*}[!t]
\centering
\includegraphics[width=\textwidth]{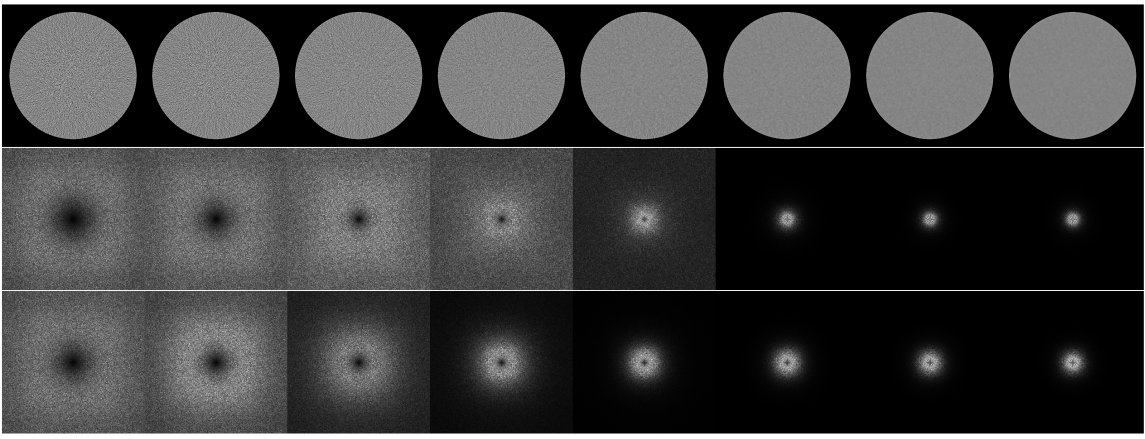}
\caption{ Simulated reconstruction paths and noise power spectrum paths of the water cylinder phantom using 100 mAs. The reconstruction path uses the PS-DOG-2, the first NPS path uses the PS-ROG-2 and the second uses the PS-DOG-2. }
\label{fig:nps}
\end{figure*}

\begin{figure*}[!t]
\centering
\includegraphics[width=0.8\textwidth]{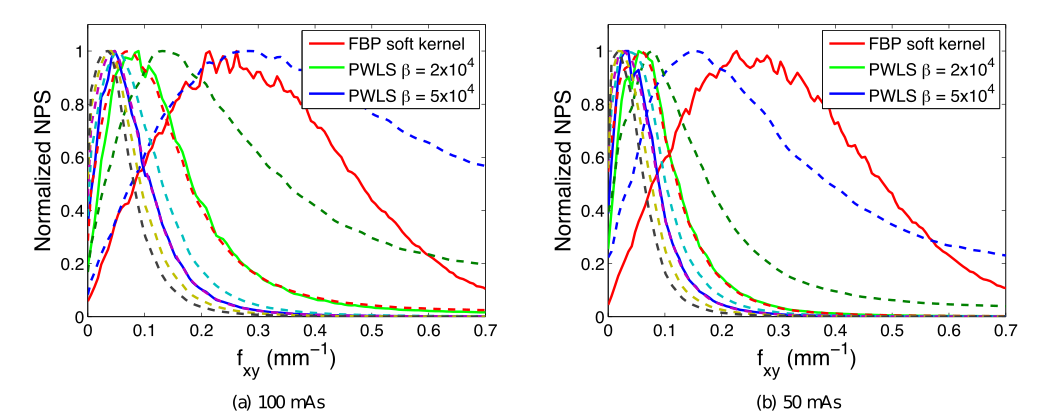}
\caption{Normalized NPS profiles of images using the filtered backprojection reconstruction, directly solved PWLS, and PS-DOG-2 (dot lines) method. }
\label{fig:nps_profile}
\end{figure*}

Figure \ref{fig:nps} shows the reconstruction paths and the corresponding in-plane noise power spectra (NPS) of the water cylinder phantom for 100 mAs scans. The peak frequency of the NPS decreases as the tuning parameter increases. A star-shaped pattern develops in the middle when the tuning parameter is big, because the penalty function in our simulations uses the differences between 6 neighboring pixels. The NPS paths are very different between the ratio-of-gradients and direction-of-gradient path seeking methods, which explains the different background textures and salt-and-pepper noise in the image domain. The ratio-of-gradients based method dose not accurately simulate the NPS of the MBIR because of the non-optimal updating set and step sizes. Figure \ref{fig:nps_profile} compares the normalized NPS profiles computed by direct optimization and the path seeking methods. For every NPS profile using the direct optimization, there is one very similar NPS profile computed by the proposed direction-of-gradient path seeking method. Note that, the 100 mAs scan has larger weighting in the data fitting function, thus the normalized NPS for 100 mAs have more high frequency components than the normalized NPS for 50 mAs with the same tuning parameter. The magnitude of noise in the 100 mAs reconstruction is not necessarily larger.

\section{Discussions}

In the MBIR problem, the number of variables is much greater than most regression problems studied in the context applied statistics. However, the methodologies of the two fields are not very different. In this paper, we proposed two path seeking algorithms that are capable of generating reasonably accurate solution path of iterative reconstruction. To the best of our knowledge, the path seeking methods are the most efficient way of computing iterative reconstructions with many different strengths of regularization. We can perform one or two reconstructions as initialization and then fill the reconstruction path without solving many large-scale optimization problems.

In the regularized linear regression such as ridge regression and LASSO \cite{Hastie2009}, regression path with one standard deviation rule provides a way of selection suitable tuning parameters. With the efficient path seeking algorithms, we suggest a path-based iterative reconstruction (PBIR) to obtain the entire reconstruction path from the same scanned data and reconstruction model. The radiologists can then select the appropriate tuning parameter in the MBIR by themselves, which is similar to the filtering kernel in the filtered backprojection reconstruction. Additionally, the path seeking algorithm allows radiologists to look at a sequence of MBIR images that can provide direct visualization of the regularization effects. 

Moreover, the PBIR can provide a fair evaluation and comparison of the reconstruction models. Over the past decade, numerous reconstruction models and regularization designs have been proposed. Different geometry discretizations, physics modeling, image roughness penalty designs, and prior knowledge augmentations provide huge variety of reconstruction models that aims to solve all kinds of challenges in the CT image reconstruction. However, the power of a model is not fully explored unless it is correctly tuned for every case and imaging task. Therefore, we suggest that one should apply the evaluation metric to the entire reconstruction path instead of to a single MBIR reconstruction. 

The PBIR also have several drawbacks. First of all, having the reconstruction path means we have to handle a much larger set of reconstruction data. The reconstruction path adds complexity in reading, storing, and transferring the images in daily practice. Secondly, the PBIR does not completely solve the problem of finding the optimal tuning parameter. To automatically select the tuning parameter, one still needs to find an appropriate image quality metric for the imaging task, which is an open research question. Although the proposed the path seeking algorithm offers good accuracy, the path images are still slightly different from the converged solution. Whether the path image is sufficient for clinical use is a question that needs to be answered. 

\section{Conclusion}

In this paper, we propose two path seeking algorithms that are capable of generating a sequence of MBIR images with different strengths of the penalty function. Simulations showed the proposed methods can produce path images that are very similar to the images computed via the direct optimization. The proposed PBIR methods enable us to obtain complete information without significant increase in the computational load. The PBIR can be easily extended to iterative image reconstructions for other image modalities such as MRI, PET, and SPECT.

\bibliographystyle{ieeetr}

\begin{thebibliography}{10}

\bibitem{Thibault2007a}
J.-B. Thibault, K.~D. Sauer, C.~A. Bouman, and J.~Hsieh, ``{A three-dimensional
  statistical approach to improved image quality for multislice helical CT},''
  {\em Med. Phys.}, vol.~34, no.~11, p.~4526, 2007.

\bibitem{Tang2009}
J.~Tang, B.~E. Nett, and G.-H. Chen, ``{Performance comparison between total
  variation (TV)-based compressed sensing and statistical iterative
  reconstruction algorithms.},'' {\em Phys. Med. Biol.}, vol.~54, no.~19,
  pp.~5781--5804, 2009.

\bibitem{Fessler:2000sl}
J.~A. Fessler, {\em {Statistical Image Reconstruction Methods for Transmission
  Tomography}}, ch.~1, pp.~1--70.
\newblock 1000 20th Street, Bellingham, WA 98227-0010 USA: SPIE, June 2000.

\bibitem{Elbakri:2002qa}
I.~A. Elbakri and J.~A. Fessler, ``{Statistical image reconstruction for
  polyenergetic X-ray computed tomography},'' {\em IEEE Trans. Med. Imaging},
  vol.~21, pp.~89--99, Feb. 2002.

\bibitem{Nuyts2013}
J.~Nuyts, B.~{De Man}, J.~A. Fessler, W.~Zbijewski, and F.~J. Beekman,
  ``{Modelling the physics in the iterative reconstruction for transmission
  computed tomography},'' {\em Phys. Med. Biol.}, vol.~58, pp.~R63--R96, June
  2013.

\bibitem{Sidky2008}
E.~Y. Sidky and X.~Pan, ``{Image reconstruction in circular cone-beam computed
  tomography by constrained, total-variation minimization.},'' {\em Phys. Med.
  Biol.}, vol.~53, no.~17, pp.~4777--4807, 2008.

\bibitem{Pfister2014}
L.~Pfister and Y.~Bresler, ``{Model-based iterative tomographic reconstruction
  with adaptive sparsifying transforms},'' in {\em Proc. SPIE Med. Imaging}
  (C.~A. Bouman and K.~D. Sauer, eds.), vol.~9020, p.~90200H, Mar. 2014.

\bibitem{Fessler:1996bs}
J.~A. Fessler and W.~L. Rogers, ``{Spatial resolution properties of
  penalized-likelihood image reconstruction: space-invariant tomographs.},''
  {\em IEEE Trans. image Process.}, vol.~5, pp.~1346--58, Jan. 1996.

\bibitem{Stayman2000}
J.~W. Stayman and J.~a. Fessler, ``{Regularization for uniform spatial
  resolution properties in penalized-likelihood image reconstruction.},'' {\em
  IEEE Trans. Med. Imaging}, vol.~19, pp.~601--15, June 2000.

\bibitem{Cho2015}
J.~H. Cho and J.~A. Fessler, ``{Regularization Designs for Uniform Spatial
  Resolution and Noise Properties in Statistical Image Reconstruction for 3-D
  X-ray CT},'' {\em IEEE Trans. Med. Imaging}, vol.~34, no.~2, pp.~678--689,
  2015.

\bibitem{chang:2015}
Z.~Chang, R.~Zhang, D.~P. {Jean-Baptiste Thibault}, L.~Fu, K.~Sauer, and
  C.~Bouman, ``{Adaptive Regularization for Uniform Noise Covariance in
  Iterative 3D CT},'' in {\em 13th Int. Meet. Fully Three-Dimensional Image
  Reconstr. Radiol. Nucl. Med.}, pp.~258--261, 2015.

\bibitem{Wang2014a}
A.~S. Wang, J.~W. Stayman, Y.~Otake, G.~Kleinszig, S.~Vogt, G.~L. Gallia, a.~J.
  Khanna, and J.~H. Siewerdsen, ``{Soft-tissue imaging with C-arm cone-beam CT
  using statistical reconstruction.},'' {\em Phys. Med. Biol.}, vol.~59, no.~4,
  pp.~1005--1026, 2014.

\bibitem{Han2015}
X.~Han, E.~Pearson, C.~Pelizzari, H.~Al-Hallaq, E.~Y. Sidky, J.~Bian, and
  X.~Pan, ``{Algorithm-enabled exploration of image-quality potential of
  cone-beam CT in image-guided radiation therapy},'' {\em Phys. Med. Biol.},
  vol.~60, no.~12, pp.~4601--4633, 2015.

\bibitem{Stayman2013}
J.~W. Stayman, H.~Dang, Y.~Ding, and J.~H. Siewerdsen, ``{PIRPLE: a
  penalized-likelihood framework for incorporation of prior images in CT
  reconstruction.},'' {\em Phys. Med. Biol.}, vol.~58, pp.~7563--7582, Oct.
  2013.

\bibitem{Wu2014a}
M.~Wu, A.~Keil, D.~Constantin, J.~Star-Lack, L.~Zhu, and R.~Fahrig, ``{Metal
  artifact correction for x-ray computed tomography using kV and selective MV
  imaging},'' {\em Med. Phys.}, vol.~41, p.~121910, 2014.

\bibitem{Dang2014b}
H.~Dang, J.~H. Siewerdsen, and J.~W. Stayman, ``{Regularization design and
  control of change admission in prior-image-based reconstruction},'' in {\em
  Proc. SPIE Med. Imaging}, vol.~9033, pp.~90330O--90330O--6, 2014.

\bibitem{Niu:2012fk}
T.~Niu and L.~Zhu, ``{Accelerated barrier optimization compressed sensing
  (ABOCS) reconstruction for cone-beam CT: phantom studies},'' {\em Med.
  Phys.}, vol.~39, pp.~4588--98, July 2012.

\bibitem{Kim2013a}
D.~Kim, S.~Ramani, and J.~A. Fessler, ``{Accelerating X-ray CT ordered subsets
  image reconstruction with Nesterov’s first-order methods},'' in {\em 12th
  Int. Meet. Fully Three-Dimensional Image Reconstr. Radiol. Nucl. Med.},
  pp.~22 -- 25, 2013.

\bibitem{Nien2014}
H.~Nien and J.~A. Fessler, ``{Fast X-Ray CT Image Reconstruction Using a
  Linearized Augmented Lagrangian Method With Ordered Subsets},'' {\em IEEE
  Trans. Med. Imaging}, vol.~34, pp.~388--399, Feb. 2015.

\bibitem{Hastie2009}
T.~Hastie, R.~Tibshirani, and J.~Friedman, {\em {The Elements of Statistical
  Learning}}.
\newblock Springer, second edi~ed., 2009.

\bibitem{Friedman2012}
J.~Friedman, ``{Fast sparse regression and classification},'' {\em Int. J.
  Forecast.}, vol.~1, 2008.

\bibitem{fiedman2009}
J.~Friedman, T.~Hastie, and R.~Tibshirani, ``{Regularization Paths for
  Generalized Linear Models via Coordinate Descent},'' {\em J. Stat. Softw.},
  vol.~30, no.~April, pp.~1--3, 2009.

\bibitem{Friedman2007pc}
J.~Friedman, T.~Hastie, H.~H\"{o}fling, and R.~Tibshirani, ``{Pathwise
  coordinate optimization},'' vol.~1, no.~2, pp.~302--332, 2007.

\bibitem{Tibshirani2011}
R.~J. Tibshirani and J.~Taylor, ``{The solution path of the generalized
  lasso},'' {\em Ann. Stat.}, vol.~39, no.~2, pp.~1335--1371, 2011.

\bibitem{Arnold2014}
T.~Arnold and R.~Tibshirani, ``{Efficient Implementations of the Generalized
  Lasso Dual Path Algorithm},'' {\em arXiv Prepr. arXiv1405.3222}, pp.~1--52,
  2014.

\bibitem{zhou2011zl}
H.~Zhou and K.~Lange, ``{A Path Algorithm for Constrained Estimation.},'' {\em
  J. Comput. Graph. Stat.}, vol.~22, pp.~261--283, Jan. 2013.

\bibitem{wu2015ps}
M.~Wu, Q.~Yang, A.~Maier, and R.~Fahrig, ``{Approximate Path Seeking for
  Statistical Iterative Reconstruction},'' in {\em Proc. SPIE Med. Imaging},
  pp.~9412--46, 2015.

\bibitem{Wu:2015tps}
M.~Wu, A.~Maier, Q.~Yang, and R.~Fahrig, ``{Improve Path Seeking Accuracy for
  Iterative Reconstruction Using the Karush-Kuhn-Tucker Conditions},'' in {\em
  Intl. Mtg. Fully 3D Image Recon. Rad.}, (New Port, RI), pp.~248 -- 251, 2015.

\bibitem{Boyd2010}
S.~Boyd and L.~Vandenberghe, {\em {Convex Optimization}}, vol.~25.
\newblock Cambridge University Press, 2010.

\bibitem{Chamboll2011}
A.~Chambolle and T.~Pock, ``{A first-order primal-dual algorithm for convex
  problems with applications to imaging},'' {\em J. Math. Imaging Vis.},
  vol.~40, no.~1, pp.~120--145, 2011.

\bibitem{hudson1994}
H.~M. Hudson and R.~S. Larkin, ``{Ordered Subsets of Projection Data},'' {\em
  IEEE Trans. Med. Imaging}, vol.~13, no.~4, pp.~601--609, 1994.

\bibitem{Kim2015}
D.~Kim, S.~Ramani, and J.~A. Fessler, ``{Combining Ordered Subsets and Momentum
  for Accelerated X-Ray CT Image Reconstruction},'' {\em IEEE Trans. Med.
  Imaging}, vol.~34, no.~1, pp.~167--178, 2015.

\bibitem{Wang2013}
A.~Wang, S.~Schafer, and J.~Stayman, ``{Soft-Tissue Imaging in Low-Dose, C-Arm
  Cone-Beam CT Using Statistical Image Reconstruction},'' in {\em Proc. SPIE
  Med. Imaging}, vol.~8668, 2013.

\end{thebibliography}

\end{document}